\documentclass[11pt,reqno]{amsart}

\usepackage{amsmath,amsthm}
\usepackage{amssymb}
\usepackage{euscript}
\usepackage{epic,eepic}

\numberwithin{equation}{subsection}

\newcommand{\sqsp}{\renewcommand{\baselinestretch}{1.2}\tiny\normalsize}

\newcommand{\relphantom[1]}{\mathrel{\phantom{#1}}}

\newtheorem{thm}[subsection]{Theorem}

\newtheorem{cor}[subsection]{Corollary}

\newtheorem*{deligne}{Deligne's Conjecture}

\theoremstyle{definition}
\newtheorem{definition}[subsection]{Definition}

\newcommand{\cat}[1]{{\EuScript #1}}

\newcommand{\cC}{\cat{C}}
\newcommand{\cH}{\cat{H}}
\newcommand{\cO}{\cat{O}}
\newcommand{\cU}{\cat{U}}
\newcommand{\cR}{\cat{R}}
\newcommand{\cT}{\cat{T}}
\newcommand{\cP}{\cat{P}}
\newcommand{\cQ}{\cat{Q}}
\newcommand{\cS}{\cat{S}}

\newcommand{\lprod}{\prec}  
\newcommand{\rprod}{\succ}  
\newcommand{\mprod}{\bot}   

\newcommand{\Cdend}{C_{\text{didend}}}

\DeclareMathOperator{\Id}{Id}
\DeclareMathOperator{\Hom}{Hom}

\begin{document}
\title[$G$-structure and Deligne's conjecture for Loday algebras]{Gerstenhaber structure and Deligne's conjecture for Loday algebras}
\author{Donald Yau}

\begin{abstract}
A method for establishing a Gerstenhaber algebra structure on the cohomology of Loday-type algebras is presented.  This method is then applied to dendriform dialgebras and three types of trialgebras introduced by Loday and Ronco.  Along the way, our results are combined with a result of McClure-Smith to prove an analogue of Deligne's conjecture for Loday algebras.
\end{abstract}

\email{dyau@math.ohio-state.edu}
\address{Department of Mathematics, The Ohio State University Newark, 1179 University Drive, Newark, OH 43055, USA}

\maketitle
\sqsp


\section{Introduction and statements of the main results}
\label{sec:intro}

The Hochschild cohomology $HH^*(A,A)$ of an associative algebra $A$ with coefficients in itself has a rich structure.  Indeed, the classical results of Gerstenhaber \cite{ger} shows that $HH^*(A,A)$ has a graded Lie bracket and a graded commutative cup product of which the Lie bracket is a graded derivation.  This structure, which is a graded version of a Poisson algebra, is now called a \emph{Gerstenhaber algebra}, or a $G$-\emph{algebra} for short.  Several other types of algebras, including coalgebras and graded associative algebras, have also been shown to admit a $G$-algebra structure in cohomology.  On the other hand, commutative algebra (Harrison) cohomology and Lie algebra (Chevalley-Eilenberg) cohomology are graded Lie algebras but not $G$-algebras \cite{gs}.  In general, whenever a new kind of algebra arises, it is an interesting and important problem to determine if it admits a $G$-algebra structure in cohomology.

Meanwhile, Loday's program \cite{loday} of studying periodicity phenomenon in algebraic $K$-theory has generated a number of new algebras.  In that program, a Lie algebra is replaced by a Leibniz algebra, which has a bracket that satisfies a version of the Jacobi identity.  If the bracket happens to be anti-symmetric, then it is a Lie algebra.  The role of an associative algebra is played by a \emph{dialgebra}, which has two associative operations satisfying three more associative-style axioms.  A dialgebra gives rise to a Leibniz algebra the same way an associative algebra gives rise to a Lie algebra.  With this analogy in mind, it is reasonable to expect that dialgebra cohomology also admits a $G$-algebra structure. The work of Majumdar and Mukherjee \cite{mm} shows that this is indeed the case.

In view of the result of \cite{mm}, it is natural to ask the question:
   \begin{quote}
   \emph{Do other Loday algebras admit a $G$-algebra structure in cohomology?}
   \end{quote}
The main purpose of this note is to extract the method used in \cite{mm}, putting it in context in such a way that it allows easy applications to other cases.  This gives a positive answer (more or less) to the question above.  This method will then be illustrated with several examples of Loday algebras.  Moreover, we combine our results with those of McClure-Smith \cite{ms} to establish a positive answer to a variation of Deligne's conjecture for Loday algebras.  We also hope that the method and examples here will make it easier to construct $G$-algebra structures on the cohomology of other Loday algebras and related algebras that may come up in the future.

A more detailed descriptions of our results follow.

\subsection{Gerstenhaber structure}
\label{subsec:main}

The method used in \cite{mm} can be summarized as follows.  The main objectives are (i) to establish a non-$\Sigma$ operad structure \cite{may1,may2} on the cochain modules and (ii) to show the existence of a \emph{multiplication} on this operad.  The operad structure arises naturally from a certain collection of functions, satisfying four conditions.  These functions are defined on certain sets that parametrize the cochain modules.  The results of Gerstenhaber and Voronov \cite{gv} about operads, brace algebras, and (homotopy) $G$-algebras  are then used to obtain the desired $G$-algebra structure on cohomology.

For many of the Loday algebras (e.g.\ (dendriform) dialgebras, associative/dendriform/cubical trialgebras \cite{lr}), there is a sequence of non-empty sets $\cU = \lbrace U_n \colon n \geq 1 \rbrace$ such that the cochain modules of an algebra $A$ of that type are given by
   \begin{equation}
   \label{eq:cochains}
   \cC(n) = C^n(A,A) = \Hom_K(K \lbrack U_n \rbrack \otimes A^{\otimes n}, A)
   \end{equation}
for $n \geq 1$, where $K$ is the ground field.  For example, in the case of dialgebras, $U_n = Y_n$ is the set of binary trees with $n+1$ leaves.  For dendriform dialgebras, $U_n = C_n = \lbrace 1, \ldots , n \rbrace$ is the finite set with $n$ elements.

\begin{definition}
\label{def:pre}
Given a sequence of non-empty sets $\cU = \lbrace U_n \colon n \geq 1 \rbrace$, define a \emph{pre-operadic system on $\cU$} to be a collection of functions
   \[
   \cR = \lbrace R_0(k; n_1, \ldots , n_k), R_i(k; n_1, \ldots , n_k) \colon k, n_1, \ldots , n_k \geq 1,\, 1 \leq i \leq k \rbrace,
   \]
where
   \[
   R_0(k; n_1, \ldots , n_k) \colon U_{n_1 + \cdots + n_k} \to U_k
   \]
and
   \[
   R_i(k; n_1, \ldots , n_k) \colon U_{n_1 + \cdots + n_k} \to U_{n_i}.
   \]
These functions are required to satisfy the following conditions:  Let $m_1, \ldots , m_N$ be positive integers, where $N = n_1 + \cdots + n_k$.  Write $N_i = n_1 + \cdots + n_i$ $(N_0 \equiv 0)$, $M_i = m_1 + \cdots + m_i$ $(M_0 \equiv 0)$, and $T_i = M_{N_i} - M_{N_{i-1}}$.  Then the functions are required to satisfy:
   \begin{description}
   \item[(1) Identity] $R_0(k; 1, \ldots , 1)$ ($k$ occurrences of $1$'s) is the identity function on $U_k$ for each $k \geq 1$.
   \item[(2) Idempotency]
      \[
      R_0(k; n_1, \ldots , n_k) R_0(N; m_1, \ldots , m_N) = R_0(k; T_1, \ldots , T_k).
      \]
   \item[(3) Commutativity] For each $i \in \lbrace 1, \ldots , k \rbrace$,
      \begin{multline*}
      R_i(k; n_1, \ldots , n_k) R_0(N; m_1, \ldots , m_N) = \\ R_0(n_i; m_{N_{i-1}+1}, \ldots , m_{N_i}) R_i(k; T_1, \ldots , T_k).
      \end{multline*}
   \item[(4) Closure] For each $i \in \lbrace 1, \ldots , k \rbrace$ and each $j \in \lbrace 1, \ldots , n_i \rbrace$,
      \[
      R_{N_{i-1}+j}(N; m_1, \ldots , m_N) = R_j(n_i; m_{N_{i-1}+1}, \ldots , m_{N_i}) R_i(k; T_1, \ldots , T_k).
      \]
   \end{description}
\end{definition}

Now let $A$ be a type of Loday algebras (e.g.\ (dendriform) dialgebras, associative/dendriform/cubical trialgebras) so that its cochain modules are given by \eqref{eq:cochains} for some sequence of non-empty sets $\cU$.  Suppose that $\cR$ is a pre-operadic system on $\cU$.  (It will be shown below that such an $\cR$ does exist for Loday algebras.)  Then for $k, n_1, \ldots , n_k \geq 1$, define maps
   \[
   \gamma \colon \cC(k) \otimes \cC(n_1) \otimes \cdots \otimes \cC(n_k) \to \cC(N)
   \]
by:
   \begin{equation}
   \label{eq:gamma def}
   \begin{split}
   \gamma(f; g_1 \otimes \cdots & \otimes g_k)(r; x_1, \ldots , x_N) \\
   & =\, f(R_0(r); g_1(R_1(r); x_1, \ldots , x_{N_1}) \otimes \cdots \\
   & \relphantom{} \relphantom{} \relphantom{} \relphantom{} \otimes g_i(R_i(r); x_{N_{i-1}+1}, \ldots , x_{N_i}) \otimes \cdots \\
   & \relphantom{} \relphantom{} \relphantom{} \relphantom{} \otimes
   g_k(R_k(r); x_{N_{k-1}+1}, \ldots , x_N)).
   \end{split}
   \end{equation}
Here $N$ and the $N_i$ are as before, $x_i \in A$, $r \in U_N$, and $R_i = R_i(k; n_1, \ldots , n_k)$ for $0 \leq i \leq k$.  Denote by $\Id_A \in \cC(1) = \Hom_K(K \lbrack U_1 \rbrack \otimes A, A)$ the canonical $1$-cochain given by
   \[
   \Id_A(r; x) = x
   \]
for $r \in U_1$ and $x \in A$.

\begin{thm}
\label{thm:main}
With the maps $\gamma$ and the $1$-cochain $\Id_A \in \cC(1)$, the collection of vector spaces $\cC = \lbrace \cC(n) \colon n \geq 1 \rbrace$ becomes a non-$\Sigma$ operad.
\end{thm}

In this case, we say that this operad is \emph{generated by the pre-operadic system $\cR$}.

A $2$-cochain $\pi \in \cC(2)$ is called a \emph{multiplication on $\cC$} if it satisfies
   \begin{equation}
   \label{eq:pi}
   \gamma(\pi; \pi, \Id_A) = \gamma(\pi; \Id_A, \pi).
   \end{equation}
Using the results and arguments of Gerstenhaber and Voronov \cite{gv}, such a multiplication generates a \emph{homotopy $G$-algebra} structure on (the brace algebra generated by) $\cC$.  Passing to cohomology, we obtain the following result:

\begin{cor}
\label{cor:main}
If $\pi \in \cC(2)$ is a multiplication on the operad $\cC$, then the corresponding cohomology
$H^*(\cC,d)$ has the structure of a $G$-algebra, where $d$ is the differential generated by $\pi$.
\end{cor}

Both Theorem \ref{thm:main} and Corollary \ref{cor:main} apply easily to the Loday algebras mentioned above.  We record it as follows.

\begin{thm}
\label{thm:examples}
Let $A$ be one of the following types of Loday algebras: dialgebras, dendriform dialgebras, associative trialgebras, dendriform trialgebras, or cubical trialgebras (so that the cochains of $A$ has the form \eqref{eq:cochains}).  Then there exists a pre-operadic system $\cR$ on $\cU$.  Moreover, the resulting operad structure on $\cC = \lbrace C^n(A,A)\rbrace$ admits a multiplication.
\end{thm}

Combining Theorem \ref{thm:main}, Corollary \ref{cor:main}, and Theorem \ref{thm:examples}, we obtain the following result in cohomology.

\begin{cor}
\label{cor:examples}
Let $A$ be a Loday algebra as in Theorem \ref{thm:examples}.  Then the cohomology $H^*(A,A)$ of $A$ has the structure of a $G$-algebra.
\end{cor}

For dialgebras, this Corollary simply recovers the results of \cite{mm}.  The other examples will be proved below in Section \ref{sec:examples}.

We now discuss a variation of Deligne's conjecture for Loday algebras.

\subsection{Deligne's conjecture for Loday algebras}
\label{subsec:Deligne}

Deligne's conjecture \cite{deligne} states that:

\begin{deligne}
The Hochschild cochain complex $C^*(A,A)$ of an associative algebra $A$ is an algebra over the singular chain operad $\cS_{\text{sing}}$ of the little squares operad $\cC_2$.
\end{deligne}

In topology, the operad $\cC_2$ is used to recognize double loop spaces and is closely related to the geometry of configuration spaces.  Deligne's conjecture, therefore, expresses a deep connection between algebra and topology.

An affirmative solution to Deligne's conjecture was given by McClure-Smith \cite{ms}, which can be summarized as follows.  There is an operad $\cH$ whose algebras are \emph{brace algebras with multiplication}, which includes $C^*(A,A)$ when $A$ is an associative algebra.  McClure and Smith showed that $\cH$ is quasi-isomorphic as a chain operad to $\cS_{\text{sing}}$.  This gives a positive answer to Deligne's conjecture.  There are also other solutions to Deligne's conjecture (see the references in \cite{ms}).

Now let $A$ be one of the types of Loday algebras in Theorem \ref{thm:examples}, so that $\cC = \lbrace C^n(A,A) \rbrace$ is an operad with multiplication.  This induces the structure of a brace algebra with multiplication on $\cC$ (Corollary \ref{cor:brace} and Theorem \ref{thm:examples}).  In particular, $\cC$ is an algebra over the operad $\cH$.  Combined with the result of McClure-Smith \cite{ms} mentioned above, this gives the following variation of Deligne's conjecture.

\begin{cor}[Deligne's conjecture for Loday algebras]
\label{cor:deligne}
Let $A$ be a Loday algebra as in Theorem \ref{thm:examples}.  Then the cochains $\cC = C^*(A,A)$ verify Deligne's conjecture: Namely, $\cC$ is an algebra over an operad $\cH$ that is quasi-isomorphic to the singular chain operad $\cS_{\text{sing}}$ of the little squares operad.
\end{cor}

It should be noted that for this Corollary to hold, we do not exactly need a Loday algebra.  More precisely, it suffices to assume that:
    \begin{enumerate}
    \item $A$ is a type of algebras whose cochain modules are in the form \eqref{eq:cochains} for some non-empty sets $U_n$.
    \item There exists a pre-operadic system on $\cU = \lbrace U_n \rbrace$.
    \item There exists a multiplication on the resulting operad structure on $\cC = \lbrace C^n(A,A) \rbrace$.
    \end{enumerate}
This might come in handy for other algebras that may come up in the future.

\subsection{Organization}

The rest of this paper is organized as follows.  The following section begins with the definition of an operad, followed by the proof of Theorem \ref{thm:main}.  It also discusses brace algebras and multiplications.  Section \ref{sec:G} discusses (homotopy) $G$-algebras, leading to Corollary \ref{cor:main}.  Our discussion on brace algebras, (homotopy) $G$-algebras, and multiplication follow Gerstenhaber and Voronov \cite{gv}.  Section \ref{sec:examples} contains a proof of Theorem \ref{thm:examples}, and hence Corollary \ref{cor:examples}, for the new examples, i.e.\ dendriform dialgebras and the three types of trialgebras.  In the final section, we show that the differential $d$ induced by the multiplication $\pi$ agrees with the differential $\delta$, up to a sign, for that particular type of algebras (Theorem \ref{thm:comparison trias}).  This ensures that the cohomology modules in Corollary \ref{cor:examples} are the intended ones.

\subsection{Acknowledgment}
The author would like to thank Mark Behrens for a discussion about this project and Jim McClure for reading an earlier version of this paper.  The author also thanks the referee for his/her helpful suggestions.


\section{Operads, brace algebras and multiplications}
\label{sec:operads}

We work over a fixed field $K$.  In this section, suppose that $A$ is a type of Loday algebras whose cochain modules are given by \eqref{eq:cochains} for some sequence of non-empty sets $\cU = \lbrace U_n \colon n \geq 1 \rbrace$.  Also, suppose that $\cR$ is a pre-operadic system on $\cU$ (see Definition \ref{def:pre}).

\subsection{Algebraic operads}
\label{subsec:operads}

A \emph{non-$\Sigma$ operad} \cite{may1,may2} is a collection $\cO = \lbrace \cO(n),\, n \geq 1 \rbrace$ of vector spaces together with structure maps
   \[
   \gamma \colon \cO(k) \otimes \cO(n_1) \otimes \cdots \otimes \cO(n_k) \,\to\, \cO(n_1 + \cdots + n_k),
   \]
for $k, n_1, \ldots , n_k \geq 1$.  The structure maps are required to satisfy the associativity condition:
   \begin{equation}
   \label{eq:operad1}
   \begin{split}
   \gamma&(\gamma(f; g_1, \ldots , g_k); h_1, \ldots , h_N) \\
   & =\, \gamma(f; \gamma(g_1; h_1, \ldots , h_{N_1}), \ldots , \\
   &\relphantom{} \relphantom{}  \relphantom{}  \relphantom{}  \gamma(g_i; h_{N_{i-1}+1}, \ldots , h_{N_i}), \ldots , \gamma(g_k; h_{N_{k-1}+1}, \ldots, h_{N_k})).
   \end{split}
   \end{equation}
Here $f \in \cO(k)$, $g_i \in \cO(n_i)$, $h_j \in \cO(m_j)$, $N = n_1 + \cdots + n_k$, and $N_i = n_1 + \cdots + n_i$.  It is also required that there be an \emph{identity element} $\Id \in \cO(1)$ such that
   \begin{equation}
   \label{eq:operad2}
   \gamma(-; \Id, \ldots , \Id) \colon \cO(k) \to \cO(k)
   \end{equation}
is the identity map.

From now on, whenever we write \emph{operad}, we mean a non-$\Sigma$ operad.

\subsection{Proof of Theorem \ref{thm:main}}
\label{subsec:proof}
Using property (1) in Definition \eqref{def:pre}, it is clear that the $1$-cochain $\Id_A$ satisfies \eqref{eq:operad2}.

For associativity \eqref{eq:operad1}, we use the notations above.  Suppose that $h_j \in \cC(m_j)$ for $1 \leq j \leq N$.  Set $N_0 = M_0 = 0$, $N_i = n_1 + \cdots + n_i$, $M_i = m_1 + \cdots + m_i$, $M = m_1 + \cdots + m_N$, and $T_i = M_{N_i} - M_{N_{i-1}}$.  Let $x_1, \ldots , x_M$ be elements of $A$, and let $r \in U_M$.  We will write $y_{s,t}$ for the sequence $y_s, \ldots , y_t$ when $s \leq t$, where $y$ can be $g$, $h$, $n$, $m$, $T$, or $x$.  Then, on the one hand, we have
   \[
   \begin{split}
   \gamma(\gamma &(f; g_{1,k}); h_{1,N})(r; x_{1,M}) \\
   & =\, f(R_0(k; n_{1,k})R_0(N; m_{1,N})(r); \cdots \\
   & \relphantom{} \relphantom{} \otimes g_i(R_i(k; n_{1,k}) R_0(N; m_{1,N})(r); \cdots \\
   & \relphantom{} \relphantom{} \relphantom{} \otimes h_{N_{i-1}+j}(R_{N_{i-1}+j}(N; m_{1,N})(r); x_{M_{N_{i-1}+j-1}+1,\, M_{N_{i-1}+j}}) \\
   & \relphantom{} \relphantom{} \relphantom{} \otimes \cdots ) \\
   & \relphantom{} \relphantom{} \otimes \cdots ).
   \end{split}
   \]
Here $1 \leq i \leq k$ and, for each $i$, $1 \leq j \leq n_i$.  In particular, not counting the $\cU$ components, the expression above displays the $i$th typical input $g_i(\cdots)$ in $f$ and the $j$th typical input $h_{N_{i-1}+j}(\cdots)$ in $g_i$.  On the other hand, we have
   \[
   \begin{split}
   \gamma &(f; \cdots , \gamma(g_i; h_{N_{i-1}+1,\, N_i}), \cdots)(r; x_{1,M}) \\
   & =\, f(R_0(k; T_{1,\, k})(r); \cdots \\
   & \relphantom{} \relphantom{} \otimes g_i(R_0(n_i; m_{N_{i-1}+1,\, N_i})R_i(k; T_{1,\, k})(r); \cdots \\
   & \relphantom{} \relphantom{} \relphantom{} \otimes h_{N_{i-1}+j}(R_j(n_i; m_{N_{i-1}+1,\, N_i}) R_i(k; T_{1,\, k})(r); x_{M_{N_{i-1}+j-1}+1,\, M_{N_{i-1}+j}}) \\
   & \relphantom{} \relphantom{} \relphantom{} \otimes \cdots ) \\
   & \relphantom{} \relphantom{} \otimes \cdots ).
   \end{split}
   \]
Comparing the first inputs (the $\cU$ components) in $f$, $g_i$, and $h_{N_{i-1}+j}$, the associativity of $\gamma$ now follows from properties (2), (3), and (4) of Definition \ref{def:pre}.

This proves that $\cC = \lbrace \cC(n) = \Hom_K(K \lbrack U_n \rbrack \otimes A^{\otimes n}, A) \colon n \geq 1 \rbrace$ is an operad. \hfill $\square$

\subsection{Brace algebras}
\label{subsec:brace}

For a graded vector space $\cO = \oplus \cO(n)$ and an element $x \in \cO(n)$, set $\deg x = n$ and $\vert x \vert = n - 1$.

Recall from \cite[Definition 1]{gv} that a \emph{brace algebra} is a graded vector space $\cO = \oplus \cO(n)$ together with a collection of braces $x \lbrace x_1, \ldots , x_n \rbrace$ of degree $-n$, satisfying
   \[
   \begin{split}
   & x\lbrace x_1, \ldots , x_m \rbrace \lbrace y_1, \ldots , y_n \rbrace \\
   & =\, \sum_{0 \leq i_1 \leq \cdots \leq i_m \leq n} (-1)^\varepsilon x \lbrace y_1, \ldots y_{i_1}, x_1 \lbrace y_{i_1 + 1}, \ldots , y_{j_1} \rbrace, y_{j_1 + 1}, \ldots , y_{i_m}, \\
   & \relphantom{} \relphantom{} \relphantom{} x_m\lbrace y_{i_m + 1}, \ldots , y_{j_m} \rbrace, y_{j_m + 1}, \ldots , y_n \rbrace.
   \end{split}
   \]
Here $\varepsilon = \sum_{p=1}^m (\vert x_p \vert \sum_{q=1}^{i_p} \vert y_q \vert)$.

According to \cite[Proposition 1]{gv}, an operad $\cC$ gives rise to a brace algebra via the braces:
   \begin{equation}
   \label{eq:operadbrace}
   x \lbrace x_1, \ldots , x_n \rbrace :=
   \sum (-1)^\varepsilon \gamma(x; \Id, \ldots , \Id, x_1, \Id, \ldots , \Id, x_n, \Id, \ldots , \Id).
   \end{equation}
Here the sum runs over all possible substitutions of $x_1, \ldots , x_n$ into $\gamma(x; \ldots )$ in the given order and $\varepsilon = \sum_{p=1}^n \vert x_p \vert i_p$, where $i_p$ is the total number of inputs in front of $x_p$.  The degree of $x \lbrace x_1, \ldots , x_n \rbrace$ is $(\sum_{p=1}^n \deg x_p) + \deg x - n$, so this operation is of degree $-n$.  This leads to the following consequence of Theorem \ref{thm:main}.

\begin{cor}
\label{cor:brace}
With the braces \eqref{eq:operadbrace}, the graded vector space $\oplus \cC(n) = \oplus C^n(A,A)$ admits the structure of a brace algebra.
\end{cor}

In such a brace algebra, define a ``comp" operation and a bracket:
   \begin{equation}
   \label{eq:circ}
   \begin{split}
   x \circ y & := x \lbrace y \rbrace, \\
   \lbrack x, y \rbrack & := x \circ y - (-1)^{\vert x \vert \vert y \vert} y \circ x.
   \end{split}
   \end{equation}
By convention, $\lbrace \rbrace$ is the identity operation, i.e., $x \lbrace \rbrace = x$.

\subsection{Multiplications}
\label{subsec:multiplications}

In an operad or a brace algebra $\cO$, a \emph{multiplication} is an element $m \in \cO(2)$ such that
   \begin{equation}
   \label{eq:mult}
   m \circ m = 0.
   \end{equation}
Since $\deg(m) = 2$, this is equivalent to
   \[
   \gamma(m; m, \Id) = \gamma(m; \Id, m).
   \]
Given such a multiplication $m$, one defines a \emph{dot product},
   \begin{equation}
   \label{eq:dot}
   x \cdot y := (-1)^{\deg x} m \lbrace x, y \rbrace,
   \end{equation}
of degree $0$ and a degree $1$ map $d$,
   \begin{equation}
   \label{eq:d}
   \begin{split}
   dx &:= \lbrack m, x \rbrack \\
      &= m \circ x - (-1)^{\vert x \vert}x \circ m.
   \end{split}
   \end{equation}
According to \cite[Proposition 2]{gv}, the map $d$ is a differential (i.e., $d^2 = 0$).  We say that $d$ is \emph{generated by $m$}.  Moreover, the dot product is associative for which $d$ is a derivation.  In particular, the dot product induces an operation on the cohomology modules defined by the differential $d$.

If $\pi \in \cC(2) = \Hom_K(K \lbrack U_2 \rbrack \otimes A^{\otimes 2}, A)$ is a multiplication, then we denote the corresponding cohomology modules by
   \[
   H^n(A,A) := H^n(\cC, d).
   \]
We will decorate this notation with a subscript for a given type of algebras.  For a given type of Loday algebras, there is usually a canonical choice of a multiplication.  The condition \eqref{eq:mult} often amounts to either the defining axioms of that type of Loday algebras or the associativity of $\pi$, as we will discuss in the examples in Section \ref{sec:examples}.

The relationships between the dot product, the bracket, and the differential are discussed next.


\section{Homotopy $G$-algebras}
\label{sec:G}

We keep the same assumptions as in the previous section.  Also, suppose that $\pi \in \cC(2) = \Hom_K(K \lbrack U_2 \rbrack \otimes A^{\otimes 2}, A)$ is a multiplication.

\subsection{Homotopy $G$-algebras}
\label{subsec:HG}

Recall from \cite[Definition 2]{gv} that a \emph{homotopy $G$-algebra} is a brace algebra $V = \oplus V^n$ with a degree $1$ differential $d$ and a degree $0$ dot product $\cdot$ that make $V$ into a differential graded associative algebra.  The dot product is required to satisfy the identity:
   \[
   (x_1 \cdot x_2) \lbrace y_1, \ldots , y_n \rbrace = \sum_{k=0}^n (-1)^\varepsilon x_1 \lbrace y_1, \ldots , y_k \rbrace \cdot x_2 \lbrace y_{k+1}, \ldots , y_n \rbrace,
   \]
where $\varepsilon = \vert x_2 \vert \sum_{p=1}^k \vert y_p \vert$.  The differential is required to satisfy the identity:
   \[
   \begin{split}
   & d(x \lbrace x_1, \ldots , x_{n+1} \rbrace) - (dx)\lbrace x_1, \ldots , x_{n+1}\rbrace \\
   & \relphantom{} - (-1)^{\vert x \vert} \sum_{i=1}^{n+1} (-1)^{\vert x_1 \vert + \cdots + \vert x_{i-1} \vert} x \lbrace x_1, \ldots , dx_i, \ldots , x_{n+1} \rbrace \\
   & =\, (-1)^{\vert x \vert \vert x_1 \vert + 1} x_1 \cdot x \lbrace x_2, \ldots , x_{n+1} \rbrace \\
   & \relphantom{} + (-1)^{\vert x \vert} \sum_{i=1}^n (-1)^{\vert x_1 \vert + \cdots + \vert x_{i-1} \vert} x \lbrace x_1, \ldots , x_i \cdot x_{i+1}, \ldots , x_{n+1} \rbrace \\
   & \relphantom - x \lbrace x_1, \ldots , x_n \rbrace \cdot x_{n+1}.
   \end{split}
   \]

A multiplication $m$ on an operad $\cO = \lbrace \cO(n) \rbrace$ gives rise to a homotopy $G$-algebra structure on the brace algebra $\oplus \cO(n)$, where the dot product and the differential are defined as in, respectively, \eqref{eq:dot} and \eqref{eq:d}.  This is Theorem 3 in \cite{gv}.  In particular, this applies to the operad $\cC$ and multiplication $\pi$.

\begin{cor}
\label{cor:HG}
With the multiplication $\pi$, the brace algebra $\oplus C^n(A,A)$ in Corollary \ref{cor:brace} admits the structure of a homotopy $G$-algebra.
\end{cor}

\subsection{$G$-algebras}
\label{subsec:G}

A \emph{$G$-algebra} \cite[2.2]{gv} is a graded vector space $H = \oplus H^n$ with a degree $0$ dot product
   \[
   - \cdot -  \colon H^m \otimes H^n \to H^{m+n}
   \]
and a degree $-1$ graded Lie bracket
   \[
   \lbrack -, - \rbrack \colon H^m \otimes H^n \to H^{m+n-1},
   \]
satisfying the following conditions:
   \begin{enumerate}
   \item The dot product is graded commutative,
   \[
   x \cdot y = (-1)^{\deg x \deg y}y \cdot x.
   \]
   \item The Lie bracket is a graded derivation for the dot product, in the sense that
   \[
   \lbrack x, y \cdot z \rbrack = \lbrack x, y \rbrack \cdot z + (-1)^{\vert x \vert \deg y} y \cdot \lbrack x, z \rbrack.
   \]
   \end{enumerate}

Corollary \ref{cor:main} in the Introduction now follows from \cite[Corollary 5]{gv} and its argument.  The dot product and the Lie bracket are induced by the ones defined in, respectively, \eqref{eq:circ} and \eqref{eq:dot}.  In particular, that the dot product is graded commutative and that the bracket is a graded derivation for the dot product are both consequences of the homotopy $G$-algebra structure in Corollary \ref{cor:HG} \cite[(8) and (9)]{gv}.


\section{Proof of Theorem \ref{thm:examples}}
\label{sec:examples}

In this section, we give a proof of Theorem \ref{thm:examples}, and hence Corollary \ref{cor:examples}.  We will not discuss dialgebras, since this was done in \cite{mm}.  In each case, we first recall some relevant definitions and the constructions of the cochain modules.

Fix a ground field $K$.

\subsection{Dendriform dialgebras}
\label{subsec:didend}

A \emph{dendriform dialgebra} $E$ over $K$ \cite[Section 5]{loday} is a $K$-vector space equipped with two binary operations,
   \begin{align*}
   \lprod \colon & E \otimes E \to E, \\
   \rprod \colon & E \otimes E \to E,
   \end{align*}
such that
   \begin{subequations}
   \label{eq:axioms}
   \begin{align}
   (x \lprod y) \lprod z &~=~ x \lprod (y \lprod z \,+\, y \rprod z), \label{eq:axiom1} \\
   (x \rprod y) \lprod z &~=~ x \rprod (y \lprod z), \label{eq:axiom2} \\
   (x \lprod y \,+\, x \rprod y) \rprod z &~=~ x \rprod (y \rprod z) \label{eq:axiom3}
   \end{align}
   \end{subequations}
for all $x, y, z \in E$.  Dendriform dialgebras are the operadic duals \cite{gk} of dialgebras.

Given a dendriform dialgebra, define a single binary operation $\ast$ by adding the two given operations:
   \[
   x \ast y \,:=\, x \lprod y + x \rprod y.
   \]
The sum of the three axioms, \eqref{eq:axiom1}, \eqref{eq:axiom2}, and \eqref{eq:axiom3}, states that $\ast$ is associative.  In particular, a dendriform dialgebra can be thought of as an associative algebra whose binary operation splits into two operations and whose associative condition splits into three conditions.

Fix a dendriform dialgebra $E$.  Denote by $C_n$ the $n$-element set $\lbrace 1, \ldots , n \rbrace$.  Define the \emph{module of $n$-cochains} as:
   \[
   \Cdend^n(E, E)
   := \Hom_K(K \lbrack C_n \rbrack \otimes E^{\otimes n}, E).
   \]
Notice that an $n$-cochain can be interpreted as an $n$-tuple of $n$-ary operations on $E$.

Using the notations in Theorem \ref{thm:main}, define
   \[
   R_0(k; n_{1,k}) \colon C_N \to C_k
   \]
by
   \[
   \label{eq:R0}
   R_0(k; n_{1,k})(r) \,=\, i \quad \text{if} \quad N_{i-1} + 1 \leq r \leq N_i
   \]
for $r \in C_N$.  For each $j \in \lbrace 1, \ldots , k \rbrace$, define
   \[
   R_j(k; n_{1,k}) \colon C_N \to C_{n_j}
   \]
by
   \[
   \label{eq:Rj}
   R_j(k; n_{1,k})(r)
   \,=\,
   \begin{cases}
   1 & \text{if } ~ 1 \leq r \leq N_{j-1}, \\
   r - N_{j-1} & \text{if } ~ 1 + N_{j-1} \leq r \leq N_j, \\
   n_j & \text{if } ~ 1 + N_j \leq r \leq N.
   \end{cases}
   \]

We claim that
   \[
   \cR = \lbrace R_0(k; n_{1,k}),\, R_j(k; n_{1,k}) \colon k, n_i \geq 1,\, 1 \leq j \leq k \rbrace
   \]
is a pre-operadic system on $\lbrace C_n \colon n \geq 1\rbrace$.  To see this, first note that $R_0(k; 1, \ldots , 1)$ ($k$ occurrences of $1$'s) is the identity function.  To prove idempotency, note that the left-hand side of condition (2) in Definition \ref{def:pre} gives
   \[
   R_0(k; n_{1,k}) R_0(N; m_{1,N})(r) \,=\, l \in \lbrace 1, \ldots , k \rbrace
   \]
if and only if
   \[
   N_{l-1} + 1 \leq R_0(N; m_{1,N})(r) \leq N_l,
   \]
which is equivalent to
   \[
   T_1 + \cdots + T_{l-1} + 1 \leq r \leq T_1 + \cdots + T_l.
   \]
This is exactly when $R_0(k; T_{1,k})(r)$ is equal to $l$, thereby proving idempotency.

With a similar reasoning, one observes that both $R_i(k; n_{1,k})R_0(N; m_{1,N})(r)$ and $R_0(n_i; m_{N_{i-1}+1,\, N_i})R_i(k; T_{1,k})(r)$ in $C_{n_i}$ are equal to:
   \[
   \begin{cases}
   1 & \text{if } ~ 1 \leq r \leq M_{N_{i-1}}, \\
   l \in \lbrace 1, \ldots , n_i \rbrace & \text{if } ~ M_{N_{i-1}+l-1} + 1 \leq r \leq M_{N_{i-1}+l}, \\
   n_i & \text{if } ~ 1 + M_{N_i} \leq r \leq M.
   \end{cases}
   \]
This proves commutativity (condition (3) in Definition \ref{def:pre}).

For closure (condition (4) in Definition \ref{def:pre}), one observes just as above that both $R_j(n_i; m_{N_{i-1}+1,N_i}) R_i(k; T_{1,k})(r)$ and $R_{N_{i-1}+j}(r)$ in $C_{m_{N_{i-1}+j}}$ are equal to:
   \[
   \begin{cases}
   1 & \text{if } ~ 1 \leq r \leq M_{N_{i-1}+j-1}, \\
   l \in \lbrace 1, \ldots , m_{N_{i-1}+j} \rbrace & \text{if } ~ r = l + M_{N_{i-1}+j-1}, \\
   m_{N_{i-1}+j} & \text{if } ~ 1 + M_{N_{i-1}+j} \leq r \leq M.
   \end{cases}
   \]
This proves closure.

Therefore, $\cR$ is a pre-operadic system on $\lbrace C_n \colon n \geq 1 \rbrace$, as claimed.  It follows from Theorem \ref{thm:main} that $\cR$ generates an operad structure on
   \[
   \cC_{\text{didend}}(E) = \lbrace \Cdend^n(E,E) \colon n \geq 1\rbrace.
   \]
Moreover, the $2$-cochain $\pi \in \Cdend^2(E,E)$ given by
   \begin{equation}
   \label{eq:pi}
   \pi(r; x \otimes y) = x \ast y = x \lprod y + x \rprod y,
   \end{equation}
for $r \in C_2$ and $x, y \in E$, is a multiplication on $\cC_{\text{didend}}(E)$.  In fact, the condition $\pi \circ \pi = 0$ is equivalent to the associativity of $\ast$.  Corollary \ref{cor:main} now implies that the corresponding cohomology,
   \[
   H_{\text{didend}}^*(E,E) = H^*(\cC_{\text{didend}}(E), d),
   \]
has a $G$-algebra structure.

\subsection{Associative trialgebras}
\label{subsec:trias}

An \emph{associative trialgebra} \cite{lr} is a vector space $A$ that comes equipped with three binary operations, $\dashv$ (left), $\vdash$ (right), and $\mprod$ (middle), satisfying the following $11$ relations for all $x, y, z \in A$:
   \begin{subequations}
   \allowdisplaybreaks
   \label{eq:triasaxioms}
   \begin{align}
   (x \dashv y) \dashv z & = x \dashv (y \dashv z), \label{axiom1} \\
   (x \dashv y) \dashv z & = x \dashv (y \vdash z), \label{axiom2} \\
   (x \vdash y) \dashv z & = x \vdash (y \dashv z), \label{axiom3} \\
   (x \dashv y) \vdash z & = x \vdash (y \vdash z), \label{axiom4} \\
   (x \vdash y) \vdash z & = x \vdash (y \vdash z), \label{axiom5} \\
   (x \dashv y) \dashv z & = x \dashv (y \,\mprod\, z), \label{axiom6} \\
   (x \,\mprod\, y) \dashv z & = x \,\mprod\, (y \dashv z), \label{axiom7} \\
   (x \dashv y) \,\mprod\, z & = x \,\mprod\, (y \vdash z), \label{axiom8} \\
   (x \vdash y) \,\mprod\, z & = x \vdash (y \,\mprod\, z), \label{axiom9} \\
   (x \,\mprod\, y) \vdash z & = x \vdash (y \vdash z), \label{axiom10} \\
   (x \,\mprod\, y) \,\mprod\, z & = x \,\mprod\, (y \,\mprod\, z). \label{axiom11}
   \end{align}
   \end{subequations}

To define the cochain modules, consider the set $\cT_n$ of planar trees with $n + 1$ leaves and one root in which each internal vertex has valence at least $2$.  We will call them \emph{trees} from now on.  The leaves of a tree $\psi \in \cT_n$ are labelled $0, 1, \ldots , n$, from left to right.  Here are the first three sets $\cT_n$:
   \[
   \begin{split}
   \cT_1 &= \lbrace \begin{picture}(16,10)     
                 \drawline(8,0)(8,2)(0,10)
                 \drawline(8,2)(16,10)
                 \end{picture}
         \rbrace, \\
   \cT_2 &= \lbrace \begin{picture}(16,10)     
                 \drawline(8,0)(8,2)(0,10)
                 \drawline(8,2)(16,10)
                 \drawline(12,6)(8,10)
                 \end{picture},\,
                 \begin{picture}(16,10)     
                 \drawline(8,0)(8,2)(0,10)
                 \drawline(4,6)(8,10)
                 \drawline(8,2)(16,10)
                 \end{picture},\,
                 \begin{picture}(16,10)     
                 \drawline(8,0)(8,10)
                 \drawline(8,2)(0,10)
                 \drawline(8,2)(16,10)
                 \end{picture}
         \rbrace, \\
   \cT_3 &= \lbrace \begin{picture}(18,11)    
                 \drawline(9,0)(9,2)(0,11)
                 \drawline(9,2)(18,11)
                 \drawline(12,5)(6,11)
                 \drawline(15,8)(12,11)
                 \end{picture},\,
                 \begin{picture}(18,11)    
                 \drawline(9,0)(9,2)(0,11)
                 \drawline(9,2)(18,11)
                 \drawline(12,5)(6,11)
                 \drawline(9,8)(12,11)
                 \end{picture},\,
                 \begin{picture}(18,11)    
                 \drawline(9,0)(9,2)(0,11)
                 \drawline(3,8)(6,11)
                 \drawline(9,2)(18,11)
                 \drawline(15,8)(12,11)
                 \end{picture},\,
                 \begin{picture}(18,11)    
                 \drawline(9,0)(9,2)(0,11)
                 \drawline(6,5)(12,11)
                 \drawline(9,8)(6,11)
                 \drawline(9,2)(18,11)
                 \end{picture},\,
                 \begin{picture}(18,11)    
                 \drawline(9,0)(9,2)(0,11)
                 \drawline(3,8)(6,11)
                 \drawline(6,5)(12,11)
                 \drawline(9,2)(18,11)
                 \end{picture},\,
                 \begin{picture}(18,11)    
                 \drawline(9,0)(9,2)(0,11)
                 \drawline(9,2)(18,11)
                 \drawline(12,5)(6,11)
                 \drawline(12,5)(12,11)
                 \end{picture},\,
                 \begin{picture}(18,11)    
                 \drawline(9,0)(9,2)(0,11)
                 \drawline(9,2)(18,11)
                 \drawline(9,2)(9,11)
                 \drawline(15,8)(12,11)
                 \end{picture},\,
                 \begin{picture}(18,11)    
                 \drawline(9,0)(9,2)(0,11)
                 \drawline(9,2)(18,11)
                 \drawline(9,2)(9,8)
                 \drawline(9,8)(6,11)
                 \drawline(9,8)(12,11)
                 \end{picture},\,
                 \begin{picture}(18,11)    
                 \drawline(9,0)(9,2)(0,11)
                 \drawline(9,2)(18,11)
                 \drawline(9,2)(9,11)
                 \drawline(3,8)(6,11)
                 \end{picture},\,
                 \begin{picture}(18,11)    
                 \drawline(9,0)(9,2)(0,11)
                 \drawline(9,2)(18,11)
                 \drawline(6,5)(6,11)
                 \drawline(6,5)(12,11)
                 \end{picture},\,
                 \begin{picture}(18,11)    
                 \drawline(9,0)(9,2)(0,11)
                 \drawline(9,2)(18,11)
                 \drawline(9,2)(6,11)
                 \drawline(9,2)(12,11)
                 \end{picture}
         \rbrace.
   \end{split}
   \]
Then the cochain modules of an associative trialgebra $A$ are defined as
   \[
   C^n_{\text{trias}}(A,A) := \Hom_K(K \lbrack \cT_n \rbrack \otimes A^{\otimes n}, A).
   \]

To define the functions $R$, first define the maps
   \begin{equation}
   \label{eq:di trias}
   d_i \colon \cT_n \to \cT_{n-1} \quad (0 \leq i \leq n),
   \end{equation}
where $d_i\psi$ is the tree obtained from $\psi$ by deleting the $i$th leaf.  These maps satisfy the simplicial relations
   \begin{equation}
   \label{eq:simplicial}
   d_i d_j = d_{j-1} d_i
   \end{equation}
for $i < j$.  Using the same kind of abbreviations and notations as in the proof of Theorem \ref{thm:main}, we define
   \[
   \begin{split}
   R_0(k; n_{1,k}) &:= d_{1,\, N_1-1} d_{N_1+1,\, N_2-1} \ldots d_{N_{k-1}+1,\, N_k-1} \colon \cT_N \to \cT_k, \\
   R_j(k; n_{1,k}) &:= d_{0,\, N_{j-1}-1} d_{N_j+1,\, N} \colon \cT_N \to \cT_{n_j}
   \end{split}
   \]
for $k, n_1, \ldots , n_k \geq 1$, $1 \leq j \leq k$.  In other words, the function $R_0$ leaves the $0$th, $N_1$th, $\ldots$ , $N_k$th leaves alone and deletes the other leaves from right to left.  The function $R_j$ leaves the $N_{j-1}$th, $(N_{j-1}+1)$st, $\ldots$ , $N_j$th leaves alone and deletes the other leaves from right to left.  These functions $R$ admit the same formulas as those in \cite[Definition 4.2]{mm}, where they are denoted by $\Gamma$ and are defined on the sets of \emph{binary trees}.

It is clear that $R_0(k; 1, \ldots , 1)$ is the identity function.  Properties (2) - (4) of Definition \ref{def:pre} are proved by the exact same argument used in \cite[Lemma 4.5]{mm}.  In fact, they all follow from the simplicial relations \eqref{eq:simplicial}.  Therefore, $\cR = \lbrace R_0(k; n_{1,k}),\, R_j(k; n_{1,k}) \colon k, n_i \geq 1,\, 1 \leq j \leq k \rbrace$ is a pre-operadic system on $\cT = \lbrace \cT_n \colon n \geq 1\rbrace$.  It follows from Theorem \ref{thm:main} that $\cC_{\text{trias}}(A) = \lbrace C^n_{\text{trias}}(A,A) \colon n \geq 1 \rbrace$ is an operad, which is generated by $\cR$.

To obtain a $G$-algebra structure, we need a multiplication.  Let $\pi \in C^2_{\text{trias}}(A,A)$ be the $2$-cochain:
   \begin{equation}
   \label{eq:pi trias}
   \pi(\psi; x, y) =
   \begin{cases}
   x \dashv y & \text{ if } \psi = \begin{picture}(16,10)     
                 \drawline(8,0)(8,2)(0,10)
                 \drawline(8,2)(16,10)
                 \drawline(12,6)(8,10)
                 \end{picture}, \\
   x \,\mprod\, y & \text{ if } \psi =  \begin{picture}(16,10)     
                 \drawline(8,0)(8,10)
                 \drawline(8,2)(0,10)
                 \drawline(8,2)(16,10)
                 \end{picture}, \\
   x \vdash y & \text{ if } \psi = \begin{picture}(16,10)     
                 \drawline(8,0)(8,2)(0,10)
                 \drawline(4,6)(8,10)
                 \drawline(8,2)(16,10)
                 \end{picture}.
   \end{cases}
   \end{equation}
Then, for $\psi \in \cT_3$, it is easy to see that the condition
   \[
   (\pi \circ \pi)(\psi; x, y, z) = 0
   \]
is equivalent to the trialgebra axioms \eqref{eq:triasaxioms}.  In fact, the $11$ possibilities of $\psi$ correspond to the $11$ trialgebra axioms.  Therefore, $\pi$ is a multiplication on the operad $\cC_{\text{trias}}(A) = \lbrace C^n_{\text{trias}}(A,A) \colon n \geq 1 \rbrace$.  It follows from Corollary \ref{cor:main} that the corresponding cohomology,
   \[
   H^*_{\text{trias}}(A,A) = H^*(\cC_{\text{trias}}(A), d),
   \]
admits a $G$-algebra structure.

\subsection{Dendriform trialgebras}
\label{subsec:tridend}

Recall from \cite{lr} that a \emph{dendriform trialgebra} is a vector space $D$ together with three binary operations, $\lprod$ (left), $\rprod$ (right), and $\cdot$ (middle), satisfying the following $7$ conditions for $x, y, z \in D$:
   \begin{subequations}
   \label{eq:tridenaxioms}
   \allowdisplaybreaks
   \begin{align}
   (x \lprod y) \lprod z &= x \lprod (y \ast z), \\
   (x \rprod y) \lprod z &= x \rprod (y \lprod z), \\
   (x \ast y) \rprod z   &= x \rprod (y \rprod z), \\
   (x \rprod y) \cdot z  &= x \rprod (y \cdot z), \\
   (x \lprod y) \cdot z  &= x \cdot (y \rprod z), \\
   (x \cdot y) \lprod z  &= x \cdot (y \lprod z), \\
   (x \cdot y) \cdot z   &= x \cdot (y \cdot z).
   \end{align}
   \end{subequations}
Here $x \ast y = x \lprod y + x \cdot y + x \rprod y$.  The operation $\ast$ is associative, which one can see by adding the seven axioms above \cite{lr}.  Dendriform trialgebras are the operadic duals of associative trialgebras.

To define the cochain modules, let $P_n$ be the set of non-empty subsets of $\lbrack n \rbrack := \lbrace 1, \ldots , n \rbrace$.  (Note: The notation in \cite{lr} is slightly different from ours.)  Then the cochain modules of $D$ are defined as
   \[
   C^n_{\text{tridend}}(A,A) := \Hom_K(K \lbrack P_n \rbrack \otimes A^{\otimes n}, A).
   \]

With the notations of the proof of Theorem \ref{thm:main}, we can define the functions $R$.  Let $X$ be an element of $P_N$.  Then
   \[
   R_0(k; n_{1,k}) \colon P_N \to P_k
   \]
is defined such that $i \in R_0(k; n_{1,k})X$ if and only if $r \in X$ for some $r$ such that $N_{i-1} + 1 \leq r \leq N_i.$  For $1 \leq j \leq k$, the function
   \[
   R_j(k; n_{1,k}) \colon P_N \to P_{n_j}
   \]
is defined by the condition: $i \in R_j(k; n_{1,k})X$ if and only if
   \[
   \begin{cases}
   r \in X \text{ for some } r \text{ such that } 1 \leq r \leq N_{j-1} + 1 & \text{ if } i = 1, \\
   i + N_{j-1} \in X & \text{ if } 2 \leq i \leq n_j - 1, \\
   r \in X \text{ for some } r \text{ such that } N_j \leq r \leq N  & \text{ if } i = n_j.
   \end{cases}
   \]

We claim that $\cR = \lbrace R_0(k; n_{1,k}),\, R_j(k; n_{1,k}) \colon k, n_i \geq 1,\, 1 \leq j \leq k \rbrace$ is a pre-operadic system on $\cP = \lbrace P_n \colon n \geq 1\rbrace$.  To see this, let $X$ be an element of $P_M$, where $M = m_1 + \cdots + m_N$.  Properties (1), (2), and (4) in Definition \ref{def:pre} follow easily by direct inspection. For property (3), suppose that $1 \leq i \leq k$ and $1 \leq l \leq n_i$.  Then one observes that $l \in R_i(k; n_{1,k})R_0(N; m_{1,N})X$ if and only if
   \[
   \begin{cases}
   r \in X \text{ for some } r \text{ with } 1 \leq r \leq M_{N_{i-1}+1} & \text{ if } l = 1, \\
   r \in X \text{ for some } r \text{ with } M_{N_{i-1}+l-1}+1 \leq r \leq M_{N_{i-1}+l} & \text{ if } 2 \leq l \leq n_i-1, \\
   r \in X \text{ for some } r \text{ with } M_{N_{i}-1}+1 \leq r \leq M & \text{ if } l = n_i.
   \end{cases}
   \]
Similarly, one observes that the above condition is equivalent to $l \in R_0(n_i; m_{N_{i-1}+1,\, N_i})R_i(k; T_{1,k})X$.  This proves property (3).  Therefore, $\cR$ is a pre-operadic system on $\cP$.  By Theorem \ref{thm:main}, $\cR$ generates an operad structure on $\cC_{\text{tridend}}(A) = \lbrace C^n_{\text{tridend}}(A,A) \colon n \geq 1 \rbrace$.

To obtain the desired $G$-algebra structure, let $\pi \in C^2_{\text{tridend}}(A,A)$ be the $2$-cochain defined by:
   \[
   \begin{split}
   \pi(X; a, b) &= a \ast b \\
   &= a \lprod b + a \cdot b + a \rprod b.
   \end{split}
   \]
The condition $\pi \circ \pi = 0$ is equivalent to the associativity of $\ast$.  Therefore, $\pi$ is a multiplication on the operad $\cC_{\text{tridend}}(A)$.  Corollary \ref{cor:main} now implies that the corresponding cohomology,
   \[
   H^*_{\text{tridend}}(A,A) = H^*(\cC_{\text{tridend}}(A), d),
   \]
admits the structure of a $G$-algebra.

\subsection{Cubical trialgebras}
\label{subsec:tricub}

Recall from \cite{lr} that a \emph{cubical trialgebra} is a vector space $A$ together with three binary operations, $\dashv$ (left), $\vdash$ (right), and $\mprod$ (middle), such that
   \begin{equation}
   \label{eq:tricubaxioms}
   (x \circ_1 y) \circ_2 z = x \circ_1 (y \circ_2 z),
   \end{equation}
for $\circ_1, \circ_2 \in \lbrace \dashv, \vdash, \mprod \rbrace$.  There are $9$ axioms in \eqref{eq:tricubaxioms}, the sum of which states that the operation $x \ast y := x \dashv y + x \,\mprod\, y + x \vdash y$ is associative.  Cubical trialgebras are operadically self-dual.

Let $Q_n$ be the set $\lbrace -1, 0, +1 \rbrace^n$.  The $i$th component of an element $X \in Q_n$ is denoted by $X_i$.  The cochain module of $A$ is defined as
   \[
   C^n_{\text{tricub}}(A,A) := \Hom_K(K \lbrack Q_n \rbrack \otimes A^{\otimes n}, A).
   \]
To define the functions $R$, let $X$ be an element of $Q_N$.  Then the function
   \[
   R_0(k; n_{1,k}) \colon Q_N \to Q_k
   \]
is given by the formula
   \[
   (R_0(k; n_{1,k})X)_i ~=~ \prod_{t=1}^{n_i} X_{N_{i-1}+t}.
   \]
for $1 \leq i \leq k$.  For $1 \leq j \leq k$, the function
   \[
   R_j(k; n_{1,k}) \colon Q_N \to Q_{n_j}
   \]
is defined by the formula
   \[
   (R_j(k; n_{1,k})X)_l = X_{N_{j-1}+l}
   \]
for $1 \leq l \leq n_j$.  We claim that $\cR = \lbrace R_0(k; n_{1,k}), R_j(k; n_{1,k}) \colon k, n_i \geq 1,\, 1 \leq j \leq k \rbrace$ is a pre-operadic system on $\cQ = \lbrace Q_n \colon n \geq 1\rbrace$.  Indeed, it is clear that $R_0(k; 1, \ldots , 1)$ is the identity function.  Properties (2) and (4) in Definition \eqref{def:pre} follow by direct inspection.  For (3), one observes that
   \[
   \begin{split}
   (R_i(k; n_{1,k})R_0(N; m_{1,N})X)_j
   &\,=\, \prod_{t=M_{N_{i-1}+j-1}}^{M_{N_{i-1}+j}} X_t \\
   &\,=\, (R_0(n_i; m_{N_{i-1}+1,\, N_i})R_i(k; T_{1,k})X)_j
   \end{split}
   \]
for $1 \leq i \leq k$ and $1 \leq j \leq n_i$.  Therefore, by Theorem \ref{thm:main}, $\cR$ generates an operad structure on $\cC_{\text{tricub}}(A) = \lbrace C^n_{\text{tricub}}(A,A) \colon n \geq 1 \rbrace$.

Let $\pi \in C^2_{\text{tricub}}(A,A)$ be the $2$-cochain given by:
   \[
   \begin{split}
   \pi(X; a, b)
   &= a \ast b \\
   &= a \dashv b + a \,\mprod\, b + a \vdash b.
   \end{split}
   \]
As in the example of dendriform trialgebras, the condition $\pi \circ \pi = 0$ is equivalent to the associativity of $\ast$.  Therefore, the corresponding cohomology,
   \[
   H^*_{\text{tricub}}(A,A) = H^*(\cC_{\text{tricub}}(A), d),
   \]
has the structure of a $G$-algebra.


\section{Comparison of cohomology}
\label{sec:comparison}

The purpose of this section is to show that the cohomology modules in Corollary \ref{cor:examples}, which arise from the differential $d$ induced by the multiplication $\pi \in \cC(2)$, are the actual ones for that particular type of algebras.  Since the various cases are rather similar, we will only work out the details in the case of associative trialgebras, which can be easily adapted to the other cases.

So let $A$ be an associative trialgebra.  The differential $\delta$ in $C^*_{\text{trias}}(A,A)$ can be figured out by considering formal deformations of associative trialgebras, along the lines of Gerstenhaber \cite{ger2}.  Deformations in the more general setting of algebras over a quadratic operad were worked out by Balavoine \cite{bal}.  We make the differential $\delta$ explicit in the associative trialgebra case to compare it with $d$.  In order to do that, we need to define a few functions on the set of planar trees.

\subsection{Functions on trees}
\label{subsec:trees}

For an element $\psi \in \cT_n$, write $\vert \psi \vert = n$.  A leaf in $\psi$ is said to be \emph{left oriented} (respectively, \emph{right oriented}) if it is the left most (respectively, right most) leaf of the vertex underneath it.  Leaves that are neither left nor right oriented are called \emph{middle leaves}.  For example, in the tree $\begin{picture}(18,11) \drawline(9,0)(9,2)(0,11) \drawline(9,2)(18,11) \drawline(9,2)(9,11) \drawline(15,8)(12,11) \end{picture}$, leaves $0$ and $2$ are left oriented, while leaf $3$ is right oriented.  Leaf $1$ is a middle leaf.

Given trees $\psi_0, \ldots , \psi_k$, their \emph{grafting} is the tree $\psi_0 \vee \cdots \vee \psi_k$ obtained by arranging $\psi_0, \ldots , \psi_k$ from left to right and joining the $k + 1$ roots to form a new (lowest) internal vertex, which is connected to a new root.  Conversely, every tree $\psi$ can be written uniquely as the grafting of $k + 1$ trees, $\psi_0 \vee \cdots \vee \psi_k$, where the valence of the lowest internal vertex of $\psi$ is $k + 1$.

For $0 \leq i \leq n+1$, define a function $\circ_i \colon \cT_{n+1} \to \lbrace \dashv, \vdash, \mprod \rbrace$
according to the following rules.  Let $\psi$ be a tree in $\cT_{n+1}$, which is written uniquely as $\psi = \psi_0 \vee \cdots \vee \psi_k$ as in the previous paragraph.  Also, write $\circ_i^\psi$ for $\circ_i(\psi)$.  Then set:
   \[
   \begin{split}
   \circ^\psi_0 & \,=\,
   \begin{cases}
   \dashv & \text{if } \vert \psi_0 \vert = 0 \text{ and } k = 1, \\
   \vdash & \text{if } \vert \psi_0 \vert > 0, \\
   \mprod & \text{if } \vert \psi_0 \vert = 0 \text{ and } k > 1,
   \end{cases} \\
   \circ^\psi_i & \,=\,
   \begin{cases}
   \dashv & \text{if the } i\text{th leaf of } \psi \text{ is left oriented}, \\
   \vdash & \text{if the } i\text{th leaf of } \psi \text{ is right oriented}, \\
   \mprod & \text{if the } i\text{th leaf of } \psi \text{ is a middle leaf},
   \end{cases} \quad (1 \leq i \leq n), \\
   \circ^\psi_{n+1} & \,=\,
   \begin{cases}
   \dashv & \text{if } \vert \psi_k \vert > 0, \\
   \vdash & \text{if } k = 1 \text{ and } \vert \psi_1 \vert = 0, \\
   \mprod & \text{if } k > 1 \text{ and } \vert \psi_k \vert = 0.
   \end{cases}
   \end{split}
   \]

\subsection{The differential $\delta$}
\label{subsec:delta}
Now the differential $\delta$ in $C^*_\text{trias}(A,A)$ is given by
   \[
   \delta^n \,=\, \sum_{i=0}^{n+1} \, (-1)^i \delta^n_i \colon C_\text{trias}^n(A,A) \,\to\, C_\text{trias}^{n+1}(A,A),
   \]
where
   \[
   (\delta^n_i f)(\psi; a_1, \ldots , a_{n+1}) \,=\,
   \begin{cases}
   a_1 \circ^\psi_0 f(d_0 \psi; a_2, \ldots , a_{n+1}) & \text{if } i = 0, \\
   f(d_i \psi; a_1, \ldots , a_i \circ^\psi_i a_{i+1}, \ldots , a_n) & \text{if } 1 \leq i \leq n, \\
   f(d_{n+1}\psi ; a_1, \ldots , a_n) \circ^\psi_{n+1} a_{n+1} & \text{if } i = n+1,
   \end{cases}
   \]
for $f \in C_\text{trias}^n(A,A)$, $\psi \in \cT_{n+1}$, and $a_1, \ldots , a_{n+1} \in A$.  Here the $d_i$ are as in \eqref{eq:di trias}.  This differential $\delta$ is similar to the one in dialgebras \cite[2.3]{fra}.

\subsection{The differential $d$}
\label{subsec:diff d}

From the construction in \S \ref{subsec:trias}, there is another differential $d$ in $C^*_\text{trias}(A,A)$ given by
   \[
   df \,=\, \pi \circ f - (-1)^{\vert f \vert} f \circ \pi,
   \]
where $\pi \in C^2_{\text{trias}}(A,A)$ is defined in \eqref{eq:pi trias}.  Now consider $d^n$ and $\delta^n \colon C^n_\text{trias}(A,A) \to C^{n+1}_\text{trias}(A,A)$.

\begin{thm}
\label{thm:comparison trias}
For each $n$, we have $d^n = (-1)^{n+1}\delta^n$.  In particular, the cohomology modules defined by $(C^*_\text{trias}(A,A), d)$ and $(C^*_\text{trias}(A,A), \delta)$ are the same.
\end{thm}

\begin{proof}
Pick an element $f \in C^n_\text{trias}(A,A)$.  Using the notations from earlier sections, we have
   \begin{equation}
   \label{eq:dnf}
   \begin{split}
   d^nf
   & \,=\, \pi \circ f + (-1)^n f \circ \pi \\
   & \,=\, \pi \lbrace f \rbrace + (-1)^n f \lbrace \pi \rbrace \\
   & \,=\, (-1)^{n-1} \gamma(\pi; \Id \otimes f) + \gamma(\pi; f \otimes \Id) \\
   & \relphantom{} \relphantom{} + \sum_{i=1}^n (-1)^{n+i-1} \gamma\left(f; \Id^{\otimes (i-1)} \otimes \pi \otimes \Id^{\otimes (n-i)}\right).
   \end{split}
   \end{equation}
We will show that these $n+2$ terms are exactly the $(-1)^i\delta^n_i$, $0 \leq i \leq n+1$, up to the sign $(-1)^{n+1}$.

Consider an element $x = \psi \otimes \mathbf{a} \in K \lbrack \cT_{n+1} \rbrack \otimes A^{\otimes (n+1)}$, where $\psi \in \cT_{n+1}$ and $\mathbf{a} = a_1 \otimes \cdots \otimes a_{n+1}$ with each $a_i \in A$.  Then
   \[
   \gamma(\pi; \Id \otimes f)(x) \,=\, \pi(R_0(2; 1, n)(\psi); a_1 \otimes f(R_2(2; 1, n)(\psi); a_{2, \, n+1})).
   \]
Using the descriptions in \S \ref{subsec:trias}, $R_2(2; 1, n)(\psi) \in \cT_n$ is obtained from $\psi$ by (1) leaving leaves $1, 2, \ldots, n+1$ alone, and (2) deleting the $0$th leaf.  That is, $R_2(2; 1, n)(\psi) = d_0 \psi$.  Likewise, the tree $R_0(2; 1, n)(\psi) \in \cT_2$ is obtained from $\psi$ by (1) leaving leaves $0, 1$, and $n+1$ alone, and (2) deleting leaves $n$, $n-1, \ldots , 2$, in this order.  Therefore, we have
   \[
   R_0(2;1, n)(\psi) \,=\,
   \begin{cases}
   \begin{picture}(16,10)     
                 \drawline(8,0)(8,2)(0,10)
                 \drawline(8,2)(16,10)
                 \drawline(12,6)(8,10)
   \end{picture}  & \text{if } \vert \psi_0 \vert = 0 \text{ and } k = 1, \\
   \begin{picture}(16,10)     
                 \drawline(8,0)(8,2)(0,10)
                 \drawline(4,6)(8,10)
                 \drawline(8,2)(16,10)
   \end{picture} & \text{if } \vert \psi_0 \vert > 0, \\
   \begin{picture}(16,10)     
                 \drawline(8,0)(8,10)
                 \drawline(8,2)(0,10)
                 \drawline(8,2)(16,10)
   \end{picture} & \text{if } \vert \psi_0 \vert = 0 \text{ and } k > 1.
   \end{cases}
   \]
This shows that
   \begin{equation}
   \label{eq:0}
   \gamma(\pi; \Id \otimes f) \,=\, \delta^n_0 f.
   \end{equation}
A similar argument shows that
   \begin{equation}
   \label{eq:n+1}
   \gamma(\pi; f \otimes \Id) \,=\, \delta^n_{n+1} f.
   \end{equation}

To finish the proof, note that for $1 \leq i \leq n$, we have
   \begin{multline*}
   \gamma\left(f; \Id^{\otimes (i-1)} \otimes \pi \otimes \Id^{\otimes (n-i)}\right)(x) \\
   \,=\, f\left(R_0(\psi); a_{1, \, i-1} \otimes \pi(R_i(\psi); a_i, a_{i+1}) \otimes a_{i+2,\, n+1}\right).
   \end{multline*}
Denoting a $k$-tuple of $1$ by $1^k$, we have $R_0(\psi) = R_0(n; 1^{i-1}, 2, 1^{n-i})(\psi)$.  This is simply $\psi$ with its $i$th leaf deleted.  That is, $R_0(\psi) = d_i \psi$.  Likewise, for $1 \leq i \leq n$, $R_i(\psi) = R_i(n; 1^{i-1}, 2, 1^{n-i})(\psi)$ is the tree in $\cT_2$ obtained from $\psi$ by (1) leaving leaves $i-1$, $i$, and $i+1$ alone, and (2) deleting the other leaves from right to left.  Therefore, we have
   \[
   R_i(n; 1^{i-1}, 2, 1^{n-i})(\psi) \,=\,
   \begin{cases}
   \begin{picture}(16,10)     
                 \drawline(8,0)(8,2)(0,10)
                 \drawline(8,2)(16,10)
                 \drawline(12,6)(8,10)
   \end{picture}  & \text{if the $i$th leaf of $\psi$ is left oriented}, \\
   \begin{picture}(16,10)     
                 \drawline(8,0)(8,2)(0,10)
                 \drawline(4,6)(8,10)
                 \drawline(8,2)(16,10)
   \end{picture} & \text{if the $i$th leaf of $\psi$ is right oriented}, \\
   \begin{picture}(16,10)     
                 \drawline(8,0)(8,10)
                 \drawline(8,2)(0,10)
                 \drawline(8,2)(16,10)
   \end{picture} & \text{if the $i$th leaf of $\psi$ is a middle leaf}.
   \end{cases}
   \]
This shows that
   \begin{equation}
   \label{eq:i}
   \gamma\left(f; \Id^{\otimes (i-1)} \otimes \pi \otimes \Id^{\otimes (n-i)}\right)
   \,=\, \delta^n_i f.
   \end{equation}
The required identity, $d^n = (-1)^{n+1} \delta^n$, is now an immediate consequence of \eqref{eq:dnf}, \eqref{eq:0}, \eqref{eq:n+1}, and \eqref{eq:i}.  This finishes the proof.
\end{proof}



\begin{thebibliography}{99}
\bibitem{bal}D.\ Balavoine, Deformations of algebras over a quadratic operad, Contemp.\ Math.\ \textbf{202} (1997), 207-234.

\bibitem{deligne}P.\ Deligne, Letter to Stasheff, Gerstenhaber, May, Schechtman, Drinfeld, May 17, 1993.

\bibitem{fra}A.\ Frabetti,  Dialgebra (co)homology with coefficients, in:  Dialgebras and related operads, 67-103, Lecture Notes in Math., \textbf{1763}, Springer, Berlin, 2001.

\bibitem{ger}M.\ Gerstenhaber, The cohomology structure of an associative ring, Ann.\ Math.\ \textbf{78} (1963), 267-288.

\bibitem{ger2}M.\ Gerstenhaber, On the deformation of rings and algebras, Ann.\ Math.\ \textbf{79} (1964), 59-103.

\bibitem{gv}M.\ Gerstenhaber and A.\ A.\ Voronov, Homotopy $G$-algebras and moduli space operad, Internat.\ Math.\ Res.\ Notices \textbf{1995}, 141-153.

\bibitem{gs}M.\ Gerstenhaber and S.\ D.\ Schack, Algebras, bialgebras, quantum groups, and algebraic deformations, in: Deformation theory and quantum groups with applications to mathematical physics (Amherst, MA, 1990), 51-92, Contemp.\ Math.\ \textbf{134}, Amer.\ Math.\ Soc., Providence, RI, 1992.

\bibitem{gk}V.\ Ginzburg and M.\ M.\ Kapranov, Koszul duality for operads, Duke Math.\ J.\ \textbf{76} (1994), 203-272.

\bibitem{loday}J.-L.\ Loday, Dialgebras, in: Dialgebras and related operads, 7-66, Lecture Notes in Math.\ \textbf{1763}, Springer, Berlin, 2001.

\bibitem{lr}J.-L.\ Loday and M.\ Ronco, Trialgebras and families of polytopes, in: Homotopy theory: relations with algebraic geometry, group cohomology, and algebraic $K$-theory, 369-398, Contemp.\ Math.\ \textbf{346}, Amer.\ Math.\ Soc., Providence, RI, 2004.

\bibitem{mm}A.\ Majumdar and G.\ Mukherjee, Dialgebra cohomology as a $G$-algebra, Trans.\ Amer.\ Math.\ Soc.\ \textbf{356} (2004), 2443-2457.

\bibitem{may1}J.\ P.\ May, The geometry of iterated loop spaces, Lectures Notes in Math.\ \textbf{271}, Springer-Verlag, Berlin-New York, 1972.

\bibitem{may2}J.\ P.\ May, Definitions: operads, algebras and modules, in:  Operads: Proceedings of Renaissance Conferences (Hartford, CT/Luminy, 1995),  1-7, Contemp.\ Math. \textbf{202}, Amer.\ Math.\ Soc., Providence, RI, 1997.

\bibitem{ms}J.\ E.\ McClure and J.\ H.\ Smith, A solution of Deligne's Hochschild cohomology conjecture, in: Recent progress in homotopy theory, 153-193, Contemp.\ Math.\  \textbf{293}, Amer.\ Math.\ Soc.\, Providence, RI, 2002.

\end{thebibliography}
\end{document}